\documentclass[11pt]{amsart}

\usepackage{amscd, amssymb}
\usepackage{enumerate}
\usepackage[all]{xypic}
\usepackage{multirow}
\usepackage{hhline}
\usepackage{verbatim}
\usepackage{mathdots}
\usepackage{comment}
\usepackage{tikz-cd}
\usepackage[titletoc]{appendix}
\usepackage{graphicx}
\usepackage{mathtools}
\usepackage{extpfeil}
\usepackage{thmtools}
\usepackage{float}

\usepackage{amsthm}
\usepackage{leftidx}
\usepackage{mathabx }
\usepackage{amsfonts}
\usepackage{bbm}

\usetikzlibrary{chains}

\allowdisplaybreaks

\makeatletter

\newcommand{\GL}{\operatorname{GL}}

\newcommand{\SO}{\operatorname{SO}}

\newcommand{\GSpin}{\operatorname{GSpin}}
\newcommand{\GPin}{\operatorname{GPin}}

\newcommand{\OO}{\operatorname{O}}

\newcommand{\la}{\langle}
\newcommand{\ra}{\rangle}

\newcommand{\GPint}{\widetilde{\operatorname{GPin}}}
\newcommand{\GSpint}{\widetilde{\operatorname{GSpin}}}

\newcommand{\OOt}{\widetilde{\operatorname{O}}}

\newcommand{\SOt}{\widetilde{\operatorname{SO}}}

\newcommand{\Fcal}{\mathcal{F}}
\newcommand{\Scal}{\mathcal{S}}
\newcommand{\Ucal}{\mathcal{U}}
\newcommand{\st}{\;:\;}

\newcommand{\sign}{\operatorname{sign}}

\theoremstyle{plain}
\newtheorem{theorem}{Theorem}[section]
\newtheorem{cor}[theorem]{Corollary}
\newtheorem{prop}[theorem]{Proposition}
\newtheorem{lemma}[theorem]{Lemma}
\newtheorem{rem}[theorem]{Remark}

\numberwithin{equation}{section}

\DeclareMathOperator{\Span}{Span}

\DeclareMathOperator{\Hom}{Hom}

\tikzset{
  symbol/.style={
    draw=none,
    every to/.append style={
      edge node={node [sloped, allow upside down, auto=false]{$#1$}}}
  }
}

\author{Melissa Emory }
\address{Department of Mathematics, Oklahoma State University; 401 Mathematical Sciences, 
Oklahoma State University
Stillwater, OK 74078}
\email{melissa.emory@okstate.edu}

\author{Yeansu Kim}
\address{Department of Mathematics Education, Chonnam National University;77 Yongbong-ro, Buk-gu, Gwangju, Korea}
\email{ykim@jnu.ac.kr}

\author{Ayan Maiti}
\address{Department of Mathematics, Purdue University; 
150 North University Street,
West Lafayette, Indiana 47907}
\email{maitia@purdue.edu}

\title[Multiplicity One for GSpin groups: Archimedean Case]{Multiplicity One Theorem for General Spin Groups:\\The Archimedean Case}

\makeatother

\begin{document}
\maketitle
\begin{abstract}
    Let $\GSpin(V)$ (resp. $\GPin(V)$) be a general spin group (resp. a general Pin group) associated with a nondegenerate quadratic space $V$ of dimension $n$ over an Archimedean local field $F$. For a nondegenerate quadratic space $W$ of dimension $n-1$ over $F$, we also consider $\GSpin(W)$ and $\GPin(W)$. We prove the multiplicity-at-most-one theorem in the Archimedean case for a pair of groups ($\GSpin(V), \GSpin(W)$) and also for a pair of groups ($\GPin(V), \GPin(W)$); namely, we prove that the restriction to $\GSpin(W)$ (resp. $\GPin(W)$) of an irreducible Casselman-Wallach representation of $\GSpin(V)$ (resp. $\GPin(V)$) is multiplicity free. 
\end{abstract}

\section{Introduction}
 Restriction problems are one of the most natural problems regarding representations, and can be formulated as follows.
 Let ${\rm G}$ be a reductive group over a local field $F$ of characteristic zero. Let ${\rm G}'$ be a reductive subgroup of ${\rm G}$, also defined over $F$. Let also $G':={\rm G}'(F)$ and $G:={\rm G}(F)$ be their $F$-points. When $\pi$ and $\pi'$ are irreducible admissible representations of $G$ and $G'$, respectively, the restriction problem asks how many times $\pi'$ appears as a quotient of $\pi$ when $\pi$ is restricted to $G'$. Formally it asks about the dimension (over $\mathbb{C}$) of the following vector space:
\[
\text{Hom}_{G'}(\pi,\pi'),
\]
where $\pi$ in the $\Hom$ space is thought of as the restriction to $G'$. There are two flavors to this problem. One is about the non-vanishing of the above Hom space, which is codified as the Gan-Gross-Prasad conjecture \cite{GGP12}. Another question concerns if the Hom space has dimension at most one, i.e. whether 
\[\text{dim}_{\mathbb{C}}(\text{Hom}_{G'}(\pi,\pi')) \leq 1\]
(which is known as a multiplicity-at-most-one or multiplicity-free theorem).
For a reductive group $G$ and its reductive subgroup $G'$, the multiplicity-at-most-one theorem in the Archimedean case (i.e., when $F$ is an Archimedean local field) explores how many times an irreducible Casselman-Wallach representation of $G'$ appears in the restriction to $G'$ of an irreducible Casselman-Wallach representation of $G$. The Archimedean case of the multiplicity-at-most-one theorem is proven for classical groups by Sun and Zhu in \cite{Sun_Zhu_Archmedean}. Using a different method of proof, Aizenbud and Gourevitch also handle the case for general linear groups in \cite{AG09b}. The main purpose of this paper is to prove the multiplicity-at-most-one theorem in the Archimedean case for a pair of two reductive non-classical groups: the general Pin groups and the general Spin groups. More specifically, let $F$ be an Archimedean local field, i.e., either $\mathbb{R}$ or $\mathbb{C}$, and let $V$ be a nondegenerate quadratic space over $F$ with dimension $n$. Correspondingly, we let $G$ be the $F$-points of either a general Pin group or a general Spin group defined over $F$, which is denoted as
\[
\GPin(V), \ \   \GSpin(V)
\]
(see Section \ref{Definition_GSpin} for exact definitions) and let $G'$ be respectively its subgroups
\[
\GPin(W), \ \ \GSpin(W)
\]
where $W \subset V$ is a nondegenerate subspace of dimension $n-1$. The multiplicity-at-most-one theorem that we prove in this paper is the following:

\begin{theorem}\label{thm:main}
Let $(G,G')$ be either $(\GPin(V),\GPin(W))$ or $(\GSpin(V),\GSpin(W))$ and let $ \pi$ be an irreducible Cassleman-Wallach representation of $G$ and $\pi'$ be an irredcuible Casselman-Wallach representation  of $G'$. Then the space of $G'$-invariant continuous bilinear functional on $\pi \times \pi'$ has dimension at most one, i.e.
$$\dim_{\mathbb{C}}\Hom_{G'}(\pi',\pi)\leq 1.$$
\end{theorem}

The technical result that will imply the above theorem (through a version of Gelfand-Kazhdan criterion) is as follows:
\begin{theorem}\label{Theorem_technical}
   Let $(G,G')$ be either $(\GPin(V),\GPin(W))$ or $(\GSpin(V),\GSpin(W))$. Then there exists a real algebraic anti-automorphism $\tau_W$ on $G'$ preserving $G$ with the following property: every generalized function on $G$  which is invariant under the adjoint action of $G'$ is automatically $\tau_W$-invariant.
\end{theorem}

The implication of Theorem 1.1  from Theorem \ref{Theorem_technical} follows from  \cite[Corollary 2.5]{SunZhu2011}. To be clear, we have 

\begin{prop}\label{prop: Thm 2 imples thm 1}
    Theorem 1.2 implies Theorem 1.1.   
\end{prop}

\begin{rem} A main ingredient in \cite{ET2023} is that they give an explicit description of the contragredient of an irreducible admissible representation of $\GSpin$ and $\GPin$ which is necessary to prove the equivalent of Proposition \ref{prop: Thm 2 imples thm 1}. Thanks to \cite[Corollary 2.5]{SunZhu2011}, this is not required in the Archimedean case.
    
\end{rem}
Let $V$ be a vector space defined over either $\mathbb{R}$ or $\mathbb{C}$ with dimension $n$. 
We denote $\widetilde{{\rm GPin}}(V):=\langle g,\beta: g\in {\rm GPin}(V) \rangle$ such that $g\beta=\beta g$ for all $g\in \GPin(V)$ and $\beta^2=1$.
Let Span$\{e_1,\dots,e_n\}=V$, where $e_{i}'$s constitute an orthogonal basis. Let
\[
W=Span\{e_1,\dots,e_{n-1}\} \text{ so that } V=W \oplus Fe \textit{ with } e:=e_n.
\]
Let $\GPint(W)\simeq\GPin(W)\rtimes\{1, e\beta\}$ act on $\GPin(V)$ (viewed merely as a set) by letting $\GPin(W)$ act by conjugation and $e\beta$ by $\tau_W$ (defined in eq 2.2). Let $\chi$ be the surjection from $\GPint(W)$ to $\{\pm 1\}$, which sends $e\beta$ to $-1$ and let $\chi$ be a continuous complex character of $\GPin(W)$. We denote $\mathcal{S}^*(\GPin(V))^{\GPint(W), \chi}$ to be the space of Schwartz distributions, in which $\GPint(W)$ acts via $\chi$. 

Theorem 1.2 can be reduced to the following vanishing assertions:

\begin{equation}\label{eq: vanishing}\mathcal{S}^*(\GPin(V))^{\GPint(W), \chi}=0;
\end{equation}
\begin{equation}\label{eq: vanishing2}\mathcal{S}^*(\GSpin(V))^{\GSpint(W), \chi}=0.
\end{equation}
Namely, every $\GPin(W)$-invariant (respectively $\GSpin$-invariant) distribution on $\GPin(V)$ (resp. $\GSpin)$ is also invariant under the involution $\tau_W$. The involution $\tau_W$ which satisfies (\ref{eq: vanishing}) and  \eqref{eq: vanishing2} is defined in (\ref{eq:tau_e_on_GPin}) and \eqref{eq:GSpininvol}.  Thanks to Proposition \ref{prop: Thm 2 imples thm 1} we only need to show the vanishing assertion of \eqref{eq: vanishing} and we are done. 

The aim of this article is to prove the above vanishing assertion in the Archimedean case. In the non-archimedean case, the multiplicity-at-most-one theorem is proven in \cite{AGRS}, \cite{Wal12} and \cite{ET2023} for the classical groups, the general Spin groups, and general Pin groups, respectively.  

This work completes the multiplicity-at-most-one theorem for the general Spin groups and the general Pin groups and completes the first step in proving the Gas-Gross-Prasad conjecture for general Spin groups.  These general Spin groups are of interest in part because they are a $\GL_1$ extension of the special orthogonal groups, and the representation theory of $\GSpin$ completely subsumes that of $\SO_n$. In addition, work on $\GSpin$ groups is needed for arithmetic application purposes as mentioned in \cite{MK2016}; orthogonal Shimura varieties are Shimura varieties of abelian type, but are finite ´etale quotients of GSpin Shimura varieties.  According to \cite{MK2016} one can easily deduce results for them from the corresponding ones for their GSpin counterparts.

The organization of the paper is as follows:
in Section \ref{Sec_BP}, we provide the background and preliminaries including  definitions of the groups $\GPin$ and $\GSpin$. In Section \ref{Sec_GPin}, we prove our main theorem for $\GPin$ groups. First, we prove that our main theorem for $\GPin$ groups reduces to Theorem \ref{thm:main_theorem_distribution}, which is the technical heart of our paper. The proof of Theorem \ref{thm:main_theorem_distribution} consists of four parts: elimination of $V$, reduction to the semisimple orbits, reduction to classical groups case, and final step of proof. Subsequently, in Section \ref{Sec_GSpin} we prove the analogous results for $\GSpin$ groups.

\subsection{Notations and conventions}

Throughout this article we assume our field $F$ to be either $\mathbb{R}$ or $\mathbb{C}$. We assume $V$ to be a vector space defined over $F$ with dimension $n$. We write $\{e_1,\dots,e_n\}$ for an orthogonal basis of $V$ and let
\[
W=Span\{e_1,\dots,e_{n-1}\} \text{ so that } V=W \oplus Fe \textit{ with } e:=e_n.
\] 
We denote $X= \GPin(V) \times V$, which will be a Nash Manifold for our purpose. Let $\frak{gpin}:= Lie(\GPin(V))$ be the Lie group of $\GPin(V)$. We denote $\widetilde{\GPin}(V):=\langle g,\beta: g\in \GPin(V) \rangle$, such that $g\beta=\beta g$ for all $g\in \GPin(V)$ and $\beta^2=1$. Therefore we have the following short exact sequence:
\[
1\longrightarrow \GPin(V)\longrightarrow\GPint(V)\xrightarrow{\;\,\chi\;\,}\{\pm1\}\longrightarrow 1,
\] 
where the surjection $\chi$ sends $\beta$ to $-1$, and $\GPint(V)\simeq\GPin(V)\times\{1, \beta\}.$ For a Nash manifold $Y$, we let $\mathcal{S}(Y)$ be the Fr\'{e}chet space of Schwartz functions on $Y$ and $\mathcal{S}^*(Y)$ be the space of Schwartz distributions on $Y$.
Denote by $\mathcal{D}(Y):=C_c^{\infty}(Y)^*$. Hence we have the following inclusion between the spaces of distributions:
\[\mathcal{S}^*(Y) \subseteq \mathcal{D}(Y).\] 
We let $P$ be the natural projection map of $\GPin$ groups onto the orthogonal groups defined in Section \ref{Definition_GSpin}. For a group $H$, we let $H_s$ be the set of semisimple elements in $H$ and we also let $H_s\slash\!\sim$ be the set of conjugacy classes in $H_s$.\\

\noindent
\textbf{Acknowledgements} This collaboration was initiated at the Midwest Representation Theory Conference at the University of Michigan in March 2022 and discussed further at the Texas-Oklahoma Representations and Automorphic Forms (TORA) Conference at University of Oklahoma in October 2023. We thank both conferences for providing a wonderful atmosphere to meet each other and collaborate.  The first author would like to thank Shuichiro Takeda and Wee teck Gan for their continued interest in the project.  The second author is grateful to Purdue University for providing excellent working conditions during a one-year research visit (July 2022 - July 2023). The second author has been supported by the National Research Foundation of Korea (NRF) grants funded by the Korea government (MSIP and MSIT) (No. RS-2022-0016551 and No. RS-2024-00415601 (G-BRL)).
\section{Background and preliminaries}\label{Sec_BP}
\subsection{Cassleman-Wallach Representation}
Let $G$ be a real reductive group and $\mathfrak{g}_{\mathbb{C}}$ be its complexified Lie algebra.
Let $(\pi, V_1)$ be a representation of $G$ and let $\mathcal{Z}(\mathfrak{g}_{\mathbb{C}})$ be the center of the universal enveloping algebra of $\mathfrak{g}_{\mathbb{C}}$. The representation $(\pi, V_1)$ is called admissible if every irreducible representation of a maximal compact subgroup $K$ of $G$ has finite multiplicity in $V_1$. The representation $(\pi, V_1)$ is called of Harish-chandra type if it is admissible and $\mathcal{Z}(\mathfrak{g}_{\mathbb{C}})$ finite.\\
The representation $(\pi, V_1)$ is called of moderate growth if the following condition holds:
for every continuous seminorm $|\cdot|_{\mu}$ on $V_1$, there exists a positive, moderate growth function $\phi$ on $G$, and a continuous seminorm $|\cdot|_{\nu}$ on $V_1$ such that
\[
|gv|_{\mu} \leq |\phi(g)|\ |v|_{\nu}, \quad \forall \ g \in G,\ v \in V_1.   
\]
The representation $(\pi, V_1)$ is called a Casselman-Wallach representation if the representation is Frechet, smooth, of moderate growth, and of Harish-Chandra type. 

\subsection{Nash groups, Nash Manifolds, and Nash maps} 
In this section we will define a Nash manifold over an Archimedean local field $F$ of characteristics zero as discussed in \cite{AAGD08}. 
\subsubsection{Semi algebraic sets}
A subset $A \subset F^n$ is called a semi-algebraic set if it satisfies the following:
there exist finitely many polynomials $f_{ij}, g_{ik} \in F[x_{1}, x_{2}, \cdots x_{n}]$ such that 
\[A = \bigcup_{i=1}^{r}\{x \in F^n: f_{i1}(x)>0, f_{i2}(x)>0, \cdots, f_{is_{i}} >0, g_{i1}(x)=0, g_{i2}(x)=0, \cdots g_{it_{i}}=0 \}.\]
In other words, the semi-algebraic sets are those that can be written as the finite union of polynomial equations and inequalities. From the definition above it is immediate that the collection of semi-algebraic sets is closed with respect to finite unions, finite intersections, and complements.
A map $\nu$ between two semi-algebraic sets $A \subset F^{n},B \subset F^m$ is called semi-algebraic if the graph of the map $\nu$ is a semialgebraic subset of $F^{m+n}$. The open semi-algebraic sets that define the topology can be realized in the following lemma:
\begin{lemma}
    Let $X \subset F^n$ be a semi-algebraic set. Then every open semi-algebraic subset of $X$ can be presented as a finite union of sets of the form $\{x \in X|f_{i}(x) > 0,i = 1...n\}$, where $f_{i}$ are polynomials in $n$ variables.
\end{lemma}
 Let $U,V \subset F^{n}$ be two open semi-algebraic sets. A smooth, semi-algebraic map between $U$ and $V$ is called a Nash map. A Nash map which is bijective, whose inverse is also a Nash map is called a Nash diffeomorphism. A Nash submanifold of $F^n$ is a semi-algebraic, smooth submanifold of $F^n$. A Nash group $G$ is a Nash manifold such that the following is a Nash map:
 \[G \times G \longrightarrow G
 \]
 \[(g,h) \mapsto gh^{-1}.\]
 Therefore we will treat the groups $\GPin(V)$ and $\GSpin(V)$ as Nash groups in the subsequent sections of this article.

\subsection{The groups $\GSpin(V)$ and $\GPin(V)$}
In this subsection we introduce the definitions of the groups $\GSpin(V)$ and $\GPin(V)$ and their properties.  A reference for the material can be found in \cite{Sch85} and \cite{Shimura}.  In the literature, the group which we refer to as $\GPin(V)$ is sometimes called the Clifford group and $\GSpin(V)$ is sometimes referred to as the special Clifford group, and are denoted by $\Gamma(V)$ and S$\Gamma(V)$, respectively.   

Here, $V$ denotes a nondegenerate quadratic space over our Archimedean local field $F$. Let $\langle -, -\rangle$ be the corresponding bi-linear form. Let $q$ be the quadratic norm defined over $F$. In this section, we also define involutions on $\GPin(V)$, which is required to prove Theorem 1.2. 

\subsubsection{Clifford Algebra} 
Let $T(V)= \bigoplus\limits_{l=0}^{\infty}V^{\otimes l}$
be the tensor algebra of $V$ and we define the Clifford algebra $C(V)$ by the following quotient:
\[
C(V)= T(V)/\langle v \otimes v - q(v)\cdot 1: v \in V\rangle.
\]
In $C(V)$ we have 
\[
v\cdot v= q(v) \quad \forall v \in V.
\]
We denote the image of $V^{\otimes l}$ in $C(V)$ as $C^{l}(V)$. Denote by
\[
C^{+}(V)= \sum_{l \text{ even}} C^{l}(V) \qquad C^{-}(V)= \sum_{l \text{ odd}} C^{l}(V)
\]
the even and odd Clifford algebras, respectively.  
Then we have the following decomposition:
\[
C(V)=C^{+}(V) \oplus C^{-}(V).
\]
The Clifford algebra is equipped with the natural involution $*$ by “reversing the indices” of $ v_{1}v_{2} \cdots v_{l} \in C^{l}(V)$, namely
\[
(v_{1}v_{2} \cdots v_{l})^{*}= v_{l}v_{l-1}\cdots v_{1}.
\]
The above involution is called the canonical involution, which preserves both $C^{+}(V)$ and $C^{-}(V)$. We define a map 

\[
\alpha : C(V) \longrightarrow C(V), \quad \alpha(x^{+}+x^{-})=x^{+}-x^{-},
\]
where $x^{+} \in C^{+}(V)$, $x^{-} \in C^{-}(V)$; in other words $\alpha$ acts on $C^{+}(V)$ as the identity and on $C^{-}(V)$ as the negative identity. For all $x \in C(V)$ the Clifford involution is 
\[
\Bar{x}=\alpha(x)^*=\alpha(x^{*}).
\]
The map sending $x$ to $\overline{x}$ is an involution on $C(V)$ and the Clifford norm is the map:
\[
N: C(V) \longrightarrow C(V), \quad N(x)=x\Bar{x}.
\]

\subsubsection{$\GPin(V)$ and $\GSpin(V)$}
\label{Definition_GSpin}
We are now in a position to define the groups $\GPin(V)$ and $\GSpin(V)$ as follows:
\[
\GPin(V):=\{g \in C(V)^{\times}: \quad \alpha(g)V g^{-1}= V\};\]
\[\GSpin(V):= \{g \in C(V)^{\times}: \quad gV g^{-1}= V\},
\]
and we call $\GSpin$ the general Spin group on $V$ and $\GPin$ the general Pin group on $V$.  We also define the projection map $P$ of $\GPin$ groups onto the orthogonal groups as follows:
\[
P(g): V \longrightarrow V, \quad P(g)v= \alpha(g)vg^{-1},
\]
for all $g \in \GPin(V)$. It is well known that $P$ surjects onto $\OO(V)$ because of the map $\alpha$. This implies that we have the following commutative diagram: 
\[
\begin{tikzcd}
  1 \arrow[r] & \GL(1) \arrow[d, equal] \arrow[r] & \GPin(V)  \arrow[r, "P"] & \OO(V) \arrow[r] & 1 \\
  1 \arrow[r] & \GL(1) \arrow[equal] \arrow[r] & \GSpin(V) \arrow[u, phantom, sloped, "\subseteq"]\arrow[r, "P"] & \SO(V)\arrow[u, phantom, sloped, "\subseteq"] \ar[r] & 1.
\end{tikzcd}
\]

\subsubsection{Involution}
Let 
\[\text{sign}: \GPin(V) \longrightarrow {\pm{1}}\]
be a homomorphism which sends the non-identity component to $-1$. Therefore the kernel of this map is $\GSpin(V)$.
We have the involution
\[
\sigma_V:\GPin(V)\longrightarrow\GPin(V),\quad
\sigma_V(g)=
\begin{cases}g^*,&\text{if $n=2k$};\\
\sign(g)^{k+1}g^*,&\text{if $n=2k-1$}.
\end{cases}
\]

Note that $\sigma_V$ preserves the semisimple conjugacy classes of $\GPin(V)$. 
We denote $\widetilde{\GPin}(V):=\langle g,\beta: g\in \GPin(V) \rangle$, such that $g\beta=\beta g$ for all $g\in \GPin(V)$ and $\beta^2=1$.
We now define the action of $\GPint(V)$ on $\GPin(V)\times V$ by
\begin{gather}\label{eq:action_on_V_GPin}
\begin{aligned}
g\cdot (h, v)&=(ghg^{-1}, P(g)v)\\
\beta\cdot (h, v)&=(\sigma_V(h), -v)
\end{aligned}
\end{gather}
for $g\in\GPin(V)$ and $(h, v)\in\GPin(V)\times V$. Note that the action of $\beta$ also preserves the semisimple conjugacy classes of $\GPin(V)$.

We assume
\[
W=Span\{e_1,\dots,e_{n-1}\} \text{ so that } V=W \oplus Fe \textit{ with } e:=e_n.
\] 
We then have $\GPin(V)_e=\GPin(W)$ and $
\GPint(V)_e=\big\la g,\, e\beta\st g\in\GPin(W)\big\ra.$
We define
\[
\GPint(W):=\GPint(V)_e.
\]
Then we also have a short exact sequence as follows:
\[
1\longrightarrow \GPin(W)\longrightarrow\GPint(W)\xrightarrow{\;\,\chi\;\,}\{\pm1\}\longrightarrow 1,
\]
where the surjection $\chi$ sends $e\beta$ to $-1$, and
$
\GPint(W)\simeq\GPin(W)\rtimes\{1, e\beta\},$
where $e\beta$ acts on $\GPin(W)$ by conjugation viewed inside $\GPint(V)$. 

We define an involution
\begin{equation}\label{eq:tau_e_on_GPin}
\tau_W:\GPin(V)\longrightarrow\GPin(V),\quad \tau_W(g)=e\sigma_V(g)e^{-1}
\end{equation}
for $g\in\GPin(V)$. This involution is precisely the action of the element $e\beta\in\GPint(V)$. Note that it is direct to show that $\tau_W(\GPin(W))=\GPin(W)$ (See \cite[Lemma 2.10]{ET2023} for more detail).

\section{$\GPin$ case}
\label{Sec_GPin}

In this section, we prove Theorem \ref{Theorem_technical} for $\GPin$ groups. Recall that Lemma \ref{mainthmtovanishing} reduces our main theorem to the following vanishing assertion:
\[
\mathcal{S}^*(\GPin(V))^{\GPint(W), \chi}=0.
\]

We first cite \cite[Theorem 2.2.5]{AG09b}, which is true for all reductive groups as follows: 
\begin{prop}
If $\Scal^*(\GPin(V))^{\widetilde{\GPin}(W), \chi}=0$,
 then $\mathcal{D}(\GPin(V))^{\widetilde{\GPin}(W), \chi}=0.$
\end{prop}
Furthermore, the following proposition is straightforward by definition:
\begin{prop}
If $\mathcal{D}(\GPin(V))^{\widetilde{\GPin}(W),\chi}=0$, then Theorem \ref{Theorem_technical} holds. Therefore the multiplicity-at-most-one theorem for $\GPin$ groups also holds.
\end{prop}

Therefore, our main theorem for $\GPin$ groups, i.e., Theorem \ref{Theorem_technical} reduces to the following theorem, which is precisely the analogue of either \cite[Theorem A]{AG09b} or \cite[Theorem 5.4]{ET2023}: 

\begin{theorem}\label{thm:main_theorem_distribution}
Let $\GPint(W)\simeq\GPin(W)\rtimes\{1, e\beta\}$ act on $\GPin(V)$ (viewed merely as a set) by letting $\GPin(W)$ act by conjugation and $e\beta$ by $\tau_W$. Then we have
\[
\mathcal{S}^*(\GPin(V))^{\GPint(W), \chi}=0.
\]
In other words, every $\GPin(W)$-invariant distribution on $\GPin(V)$ is also invariant under the involution $\tau_W$.
\end{theorem}

The rest of this section is to prove
Theorem \ref{thm:main_theorem_distribution}, which is the technical heart of the paper. We adapt the arguments in \cite[Section 7]{ET2023} to our case and it consists of three steps of reductions.

\subsection{Reduction I: Elimination of $W$}
In this subsection, We reduce Theorem \ref{thm:main_theorem_distribution} to the follwoing vanishing assertion:
\begin{equation}\label{eq:chi_distribution_zero}
\Scal^*(\GPin(V)\times V)^{\GPint(V),\chi}=0.
\end{equation}
The key ingredient is the following version of Frobenius descent \cite[Theorem 2.5.7]{AG09a}.

\begin{lemma}[Frobenius descent]\label{lemma:Frobenius_descent}
Let $G$ be a Nash group which is unimodular. Let $X$ and $Y$ be Nash manifolds on which $G$ acts. Further assume that the action of $G$ on $Y$ is transitive. Suppose we have a continuous $G$-equivariant Nash map
\[
\phi:X \rightarrow Y,
\]
namely $\varphi(g\cdot x)=g\cdot\varphi(x)$ for all $g\in G$ and $x\in X$. Fix $y\in Y$. Assume the stabilizer $G_y\subseteq G$ of $y$ is unimodular which implies that there exists a $G$-invariant measure on $Y$.  Fix this measure.  Let $\chi:G\to\mathbb{C}^1$ be a character of $G$.  Then there is a canonical isomorphism
\[
\mathcal{S}^*(X)^{G,\chi}\simeq \mathcal{S}^*(\phi^{-1}(y))^{G_y,\chi}.
\]
\end{lemma}

We are now ready to prove the following:
\begin{prop}\label{prop:removing_W_from_vanishing_assertion}
If $\mathcal{S}^*(\GPin(V)\times V)^{\GPint(V), \chi}=0$ then $\mathcal{S}^*(\GPin(V))^{\GPint(W)), \chi}=0$.
\end{prop}
\begin{proof}
We follow the main argument of the proof in \cite[Proposition 7.2]{ET2023} and adapt those to Archimedean version. Briefly, we write the main ideas of the proof but, for completeness, we write in detail in case the proof is not certain in the Archimedean case. 

Our goal is to prove 
\[
\mathcal{S}^*(\GPin(V))^{\GPint(W), \chi}\subseteq \Scal^*(\GPin(V)\times V)^{\GPint(V), \chi}.
\]

Let $X:=\{(g,v) \in \GPin(V)\times V \st\langle v,v\rangle=\langle e,e \rangle \}.$ One can see that $X$ is invariant under $\GPint(V)$ and we have the following:
\begin{lemma}
    $X$ is closed in $\GPin(V)\times V$.
\end{lemma}
\begin{proof}
    We have the quadratic form $q :V\rightarrow F$ and we know that $q(v)=v^2$ and $v^2 \in F$.  Hence $\langle v,v \rangle = \langle e,e\rangle $ implies 
\[
q(v+v)-q(v)-q(v)=q(e+e)-q(e)-q(e),
\]
\[4v^2-2v^2=4e^2-2e^2,
\]
\[v^2=e^2.
\]
Since $v^2,e^2 \in F$, this is a discrete set.  Therefore, $X$ is closed in $\GPin(V) \times V$.
\end{proof}

Therefore, we have
\[
\mathcal{S}^*(X)^{\GPint(V), \chi}\subseteq\mathcal{S}^*(\GPin(V)\times V)^{\GPint(V), \chi}.
\]

Let $Y:=\{v\in V\st \langle v,v\rangle=\langle e,e \rangle\}.$ One can also see that $Y$ is invariant under the action of $\GPin(V)$. By Witt's theorem, we know $\OO(V)$ acts transitively on $Y$ and hence $\GPint(V)$ acts transitively on $Y$.

Now, consider the projection
\[
\phi:X\longrightarrow Y,\quad (g, v)\mapsto v,
\]
which is $\GPint(V)$-equivariant. Recall that the stabilizer $\GPint(V)_e=\GPint(W)$ of $e$ is unimodular. Hence by the Frobenius descent (Lemma \ref{lemma:Frobenius_descent}) applied to this $\phi$, we obtain the canonical isomorphism
\begin{equation}\label{eq:canonical_iso_Forb1}
\mathcal{S}^*(X)^{\GPint(V), \chi} \simeq \mathcal{S}^*(\phi^{-1}(e))^{\GPint(V)_e,\chi}.
\end{equation}
By the obvious identification $\GPin(V)\times\{e\}\simeq\GPin(V)$ of sets, we have
\[
\mathcal{S}^*(\phi^{-1}(e))^{\GPint(V)_e,\chi}\simeq\mathcal{S}^*(\GPin(V))^{\GPint(W),\chi}.
\]
Hence we have
\[
\mathcal{S}^*(\GPin(V))^{\GPint(W),\chi}\simeq \mathcal{S}^*(X)^{\GPint(V), \chi}\subseteq \mathcal{S}^*(\GPin(V)\times V)^{\GPint(V), \chi}.
\]
The proposition is proven.
\end{proof}
\subsection{Reduction II: the semisimple orbits}
The proof of our main theorem is now reduced to showing the vanishing assertion (\ref{eq:chi_distribution_zero}). 
In this subsection we further reduce the vanishing assertion to the classical group scenarios as presented in \cite{Sun_Zhu_Archmedean} or in \cite{Wal12}.

The main idea is that any distribution in $\mathcal{S}^*({\GPin(V)}\times V)^{\widetilde{\GPin(V)},\chi}$ is supported in a smaller set through Harish-Chandra's descent and Bernstein's localization principle. We cite those theorems for the Archimedean case as follows:

\begin{theorem}\label{localization} 
    (Localization principle). Let a real reductive group $G$ act on a smooth algebraic variety $X$. Let $Y$ be an algebraic variety and $\phi: X \rightarrow Y$ be an affine algebraic $G$-invariant map. Let $\chi$ be a character of $G$. Suppose that for any $y \in Y(F)$ we have $\mathcal{D}_{X(F)}\left(\phi^{-1}(y)(F)\right)^{G, \chi}=0$. Then $\mathcal{D}(X(F))^{G, \chi}=0$.    
\end{theorem}
\begin{proof}
See \cite[Corollary A.0.1]{AG09b}.
\end{proof}
\begin{rem} \label{rem: D and S are interchangeable} We can swiftly interchange the notation $\mathcal{D}$ and $\mathcal{S}^*$, due to the proof of \cite[Theorem 4.0.2]{AG09b}.
\end{rem}
\begin{cor}\label{localizationcor}
     Let a real reductive group $G$ act on a smooth algebraic variety $X$. Let $Y$ be an algebraic variety and $\phi: X \rightarrow Y$ be an affine algebraic $G$-invariant submersion. Suppose that for any $y \in Y(F)$ we have $\mathcal{S}^*\left(\phi^{-1}(y)\right)^{G, \chi}=0$. Then $\mathcal{D}(X(F))^{G, \chi}=0$. 

\end{cor}
\begin{proof}
    See \cite[Corollary A.0.3]{AG09b}.
\end{proof}

\begin{prop}\label{prop:reduction_to_semisimple_orbit}

Define a map
\[
\theta:\GPin(V)\times V\longrightarrow \OO(V)_{s}\slash\!\sim
\]
by
\[
(g, v)\mapsto P(g)_s\slash\!\sim,
\]
where $P(g)_s\slash\!\sim$ is the conjugacy class of the semisimple part of $P(g)$ under the Jordan decomposition.

If
\[
\Scal^*(\theta^{-1}(\gamma))^{\GPint(V), \chi}=0
\]
for each semisimple conjugacy class $\gamma \in \OO(V)_{s}\slash\!\sim$, then
\[
\mathcal{S}^*(\GPin(V)\times V)^{\GPint(V),\chi}=0.
\]
\end{prop}
\begin{proof}
We follow the main argument of the proof in \cite[Proposition 7.4]{ET2023} and we briefly write the main ideas of the proof but in case the proof is different we write in the  details. 
Let $Y$ be the space of polynomials of degree at most $n=\dim V$, which is a smooth algebraic variety and let $
\phi:\GPin(V)\times V \longrightarrow Y$ be as in proof of \cite[Proposition 7.4]{ET2023}. One can see that $\phi$ is an affine algebraic $G$-invariant submersion. Let $f\in Y$ be a polynomial. Since the fiber $\phi^{-1}(f)$ is preserved by $\GPint(V)$, Theorem \ref{localization} implies that if $\Scal^*(\phi^{-1}(f))^{\GPint(V),\chi}=0$
for all $f\in Y$, then $\Scal^*(\GPin(V)\times V)^{\GPint(V),\chi}=0$. Furthermore, we have $\phi^{-1}(f)=\Fcal_f\times V$
where $
\Fcal_f=\{g\in\GPin(V)\st \text{the char.\ poly.\ of $P(g)$ is $f$}\},$ which is a Nash manifold.

To use the Bernstein localization principle, we need to show that each element in $P(\Fcal_f)_{s}\slash\!\sim$ is closed and that $P(\mathcal{F}_{f})_{s}/\sim$ is an embedded submanifold of  $P(\mathcal{F}_{f})_{s}$ (basically it means the topology is a subspace topology).  The first claim follows from \cite[Proposition 10.1]{BH62}. The second one is well-known and can be followed from \cite[Theorem 2.13]{MZ40}.

We can now consider the map $
\theta:\Fcal_f\times V\longrightarrow P(\Fcal_f)\longrightarrow P(\Fcal_f)_{s}\longrightarrow P(\Fcal_f)_{s}\slash\!\sim,
$ exactly as in the proof of \cite[Proposition 7.4]{ET2023}. Then this map $\theta$ is indeed an affine algebraic map. In the Archimedean case, the involution $\sigma_V$ also preserves the semisimple conjugacy classes. Therefore, for each semisimple conjucagy class $\gamma \in P(\Fcal_f)_{s}\slash\!\sim$ the fiber $\theta^{-1}(\gamma)$ is invariant under $\GPint(V)$ and Theorem \ref{localization} implies that if $\Scal^*(\theta^{-1}(\gamma))^{\GPint(V),\chi}=0$ for all semisimple conjugacy class $\gamma$ of $\OO(V)$, then $\Scal^*(\GPin(V)\times V)^{\GPint(V),\chi}=0,$
which completes the proof of the Proposition.

\end{proof}

\subsection{Reduction III: $\OO(V)$ situation}
In this section we further reduce the vanishing assertion of the hypothesis in Proposition \ref{prop:reduction_to_semisimple_orbit} further to the orthogonal groups situation \cite{Sun_Zhu_Archmedean}.

\begin{rem}
  Let $\Ucal\subseteq\GPin(V)$ be the set of unipotent elements in $\GPin(V)$. Then, for each $g\in\GPin(V)$ both $\Ucal$ and $\Ucal_g$ are closed as smooth algebraic variety in $\GPin(V)$. It is also well known that the restriction to $\Ucal$ of the canonical projection $P:\GPin(V)\to\OO(V)$ is one-to-one, which allows us to identify the set of unipotent elements in $\OO(V)$ with $\Ucal$.  
\end{rem}

For a semisimple conjugacy class $\gamma\subseteq\OO(V)_s$ of $P(g)$, we cosnider the following map exactly as in \cite{ET2023}:
\[
\theta:\theta^{-1}(\gamma)\longrightarrow \gamma, \quad (g, v)\mapsto P(g)_s,
\]
Applying the Frobenius descent to the above map, we have 
\[
\Scal^{*}(\theta^{-1}(\gamma))^{\GPint(V), \chi}\simeq\Scal^{*}(Z^\circ g\,\Ucal_g\times V)^{\GPint(V)_{g}, \chi}
\]
since $\theta^{-1}(P(g))=Z^\circ g\,\Ucal_g\times V$.

Furthermore, applying the Bernstein localization principle, we have the following:

\begin{lemma}
If $\Scal^{*}(g\,\Ucal_g\times V)^{\GPint(V)_{g}, \chi}=0$
for all semisimple $g\in\GPin(V)$,
then $\Scal^{*}(Z^\circ g\,\Ucal_g\times V)^{\GPint(V)_{g}, \chi}=0$.
\end{lemma}
\begin{proof}
First, since $zg$ is semisimple for all $z\in Z^\circ$ and $\Ucal_{zg}=\Ucal_g$, we have 
\[
\Scal^{*}(zg\,\Ucal_g\times V)^{\GPint(V)_{g}, \chi}=0
\]
for all $z\in Z^\circ$.

We then consider the following map:
\[
Z^\circ g\,\Ucal_g\times V\longrightarrow Z^\circ,\quad zg\, \Ucal_g\mapsto z.
\]
Now, as each fiber $zg\,\Ucal_g\times V$ is preserved by $\GPint(V)_g$, we may apply the Bernstein Localization principle (Corollary \ref{localizationcor}) to arrive at the desired conclusion. 
\end{proof}

We finally reduce our main theorem to the $\OO(V)$ situation. Recall that $\GPin(V)$ and $\OO(V)$ have the same set of unipotent elements. This implies that for each semisimple $g$ we have the bijection
\[
g\,\Ucal_g\longrightarrow P(g\,\Ucal_g)
\]
induced by the canonical projection $P$. Furthermore, this map intertwines the actions of $\GPint(V)_g$ and $P(\GPint(V)_g)$. Note that the kernel $Z_0$ of $P$ act trivially on the space $\Scal^{*}(g\,\Ucal_g\times V)$. Therefore, we have 
\[
\Scal^{*}(g\,\Ucal_g\times V)^{\GPint(V)_{g}, \chi}\simeq\Scal^{*}(P(g\,\Ucal_g)\times V)^{P(\GPint(V)_g), \chi}.
\]
Hence to show our main theorem, it suffices to show the following $\OO(V)$ situation \cite[Theorem A]{Sun_Zhu_Archmedean}:
\begin{equation}\label{eq:vanishing_for_orthogonal}
\Scal^{*}(P(g\,\Ucal_g)\times V)^{P(\GPint(V)_g), \chi} \subseteq \Scal^{*}(P(\GPin(V)_g)\times V)^{P(\GPint(V)_g), \chi}=0
\end{equation}
for all semisimple $g\in\GPin(V)$. 

\subsection{End of proof}

In case $P(\GPin(V)_g)=\OO(V)_{P(g)}$, the vanishing assertion \eqref{eq:vanishing_for_orthogonal} is equivalent to show $\Scal^{*}(\OO(V)_{P(g)}\times V)^{\OO(V)_{P(g)}, \chi}=0$, which is precisely the assertion proven in \cite{Sun_Zhu_Archmedean}. However, we do not always have $P(\GPin(V)_g)=\OO(V)_{P(g)}$ as shown in \cite[Lemma 3.11]{ET2023}. Accordingly, we need to modify \cite{Sun_Zhu_Archmedean}. The difference is that $P(\GPin(V)_g)$ might have a factor of $\SO$ as in \cite[Lemma 3.11]{ET2023}, for which we need the result of Sun-Zhu \cite{Sun_Zhu_Archmedean} for the $\SO$ case. The complete argument in the case of locally compact totally disconnected space is present in \cite{ET2023}. The same argument works in our case and we do not repeat the same argument. 

In conclusion, we have shown the vanishing assertion of \eqref{eq: vanishing} which gives Theorem \ref{thm:main} for the groups $\GPin (V)$ and $\GPin (W)$.

\quad

\section{GSpin case}
\label{Sec_GSpin}
In this section, we prove Theorem \ref{Theorem_technical} for $\GSpin(V)$. The proof follows the same line as the GPin case but we need to make appropriate changes. Note that we follow the main arguments in \cite[Proposition 9]{ET2023} and we adapt its arguments to Archimedean case and we briefly explain the main ideas. To point out the difference, you can see the proof of Proposition \ref{vanishing_inc_GSpin}. In particular, we define a group $\GSpint(V)$ which is the analogue of $\SOt(V)$ following exactly as the non-archimedean case.  
First, we recall our basic setup. As before, $V$ is a quadratic space with
\[
\dim_FV=n=\begin{cases}2k;\\ 2k-1.\end{cases}
\]
We fix an orthogonal basis $e_1,\dots,e_{n-1}, e_n$, and assume $W=\Span\{e_1,\dots,e_{n-1}\}$ so that $V=W\bigoplus Fe$ with $e:=e_n$. 
The group $\GSpint(V)$ we define as

\[
\GSpint(V)=\big\la g,\, e^k\beta\st g\in\GSpin(V)\big\ra\subseteq\GPint(V),
\]
so that we have
\[
1\longrightarrow \GSpin(V)\longrightarrow\GSpint(V)\xrightarrow{\;\,\chi\;\,}\{\pm1\}\longrightarrow 1.
\]
The surjection $\chi$ sends $e^k\beta$ to $-1$, and
\[
\GSpint(V)\simeq\GSpin(V)\times\{1, e^k\beta\},
\]
where $e^k\beta$ acts on $\GSpin(V)$ by conjugation viewed inside $\GPint(V)$. Since $\GSpint(V)$ is a subgroup of $\GPint(V)$, it acts on $\GSpin(V)\times V$ (viewed merely as a set) by restricting the action of $\GPint(V)$ as
\begin{gather}\label{eq:action_of_GSpint}
\begin{aligned}
g\cdot (h, v)&=(ghg^{-1}, P(g)v)\\
e^k\beta\cdot (h, v)&=(e^k\sigma_V(h)e^{-k}, -P(e)^kv),
\end{aligned}
\end{gather}
where $(h, v)\in\GSpin(V)\times V$.

We let $\GSpint(V)_e$ be the stabilizer of $e\in V$ under the action of $\GSpint(V)$ on $V$ as usual. Analogously to the $\SO(V)$ case, one can then show
\[
\GSpint(V)_e=\big\la g, e_{n-1}^{k-1}e\beta\st g\in\GSpin(W)\big\ra,
\]
and we define
\[
\GSpint(W):=\GSpint(V)_e.
\]
We have
\[
1\longrightarrow \GSpin(W)\longrightarrow\GSpint(W)\xrightarrow{\;\,\chi\;\,}\{\pm1\}\longrightarrow 1,
\]
where the surjection $\chi$ sends $e_{n-1}^{k-1}e\beta$ to $-1$, and
\[
\GSpint(W)\simeq\GSpin(W)\rtimes\{1, e_{n-1}^{k-1}e\beta\},
\]
where the action of $e_{n-1}^{k-1}e\beta$ is by conjugation viewed inside $\GPint(V)$.

We define an involution 
\begin{equation}\label{eq:GSpininvol}
\tau_W:\GSpin(V)\longrightarrow\GSpin(V),\quad \tau_W(g)=(e_{n-1}^{k-1}e)\sigma_V(g)(e_{n-1}^{k-1}e)^{-1},
\end{equation}
for $g\in\GSpin(V)$. This is the action of $e_{n-1}^{k-1}e\beta\in\GSpint(V)_e$ on $\GSpin(V)$. Since $e$ commutes with all the elements in $\GSpin(W)$, we have $\tau_W(\GSpin(W))=\GSpin(W)$.

We have the canonical projection
\[
P:\GSpint(V)\longrightarrow \SOt(V),\quad g\mapsto P(g),\;e^k\beta\mapsto r_e^k\beta,
\]
which is nothing but the restriction of the canonical projection $P:\GPint(V)\to\OOt(V)$. We then have
\[
P(\GSpint(V)_e)=\SOt(V)_e.
\]

Let $g\in\GSpin(V)$ be semisimple, and set $h:=P(g)\in\SO(V)$. If
\[
\OO(V)_h\simeq G_1\times\cdots\times G_m\times\OO(V_+)\times\OO(V_-)
\]
as before, then
\[
\SO(V)_h\simeq G_1\times\cdots\times G_m\times S(\OO(V_+)\times\OO(V_-)),
\]
where
\[
S(\OO(V_+)\times\OO(V_-))=(\OO(V_+)\times\OO(V_-))\cap\SO(V_+\oplus V_-)
\]
by \cite[Proposition A.4]{ET2023}. Following exactly as in the proof of \cite[Lemma 9.1]{ET2023}, we have the following Lemma:
\begin{lemma}\label{lemma:centralizer_GSpin}
Keeping the above notation, we have
\[
P(\GSpin(V)_g)\simeq G_1\times\cdots\times G_m\times \SO(V_+)\times\SO(V_-) \subset \SO(V)_h.
\]
\end{lemma}

\subsection{Vanishing of distribution}
Analogously to the GPin case, we prove the following main technical result:
\begin{equation}\label{eq:vanishing_distribution_GSpin}
\Scal^*( \GSpin(V))^{\GSpint(W), \chi}=0,
\end{equation}
where $\GSpint(W)\simeq \GSpin(W)\times\{1, e_{n-1}^{k-1}e\beta\}$ acts on $\GPin(V)$ by restricting the actions \eqref{eq:action_of_GSpint}. In particular, the element $e_{n-1}^{k-1}e\beta$ acts via the involution $\tau_W$, which preserves $\GSpin(W)$ setwise. Recall the action of $\GSpint(V)$ on $\GSpin(V)\times V$ is defined in \eqref{eq:action_of_GSpint}.

\begin{prop}\label{vanishing_inc_GSpin}
We have a natural inclusion
\[
\Scal^*(\GSpin(V))^{\GSpint(W), \chi}\subseteq \Scal^*(\GSpin(V)\times V)^{\GSpint(V), \chi}.
\]
Hence if
\[
\Scal^*(\GSpin(V)\times V)^{\GSpint(V), \chi}=0,
\]
then $\Scal^*(\GSpin(V))^{\GSpint(W), \chi}=0$.
\end{prop}
\begin{proof}
This can be proven in the same way as Proposition \ref{prop:removing_W_from_vanishing_assertion}. Namely let
\begin{align*}
X&:=\{(g,v)\in\GSpin(V)\times V\st \la v, v\ra=\la e, e\ra\}\\
Y&:=\{v\in V\st \la v, v\ra=\la e, e\ra\},
\end{align*}
and consider the projection
\[
\phi:X\longrightarrow Y.
\]
By Witt's theorem, $\GSpin(V)$ acts transitively on $Y$.  Moreover, using the actions listed in (\ref{eq:action_of_GSpint}), one can show that $\phi(g \cdot (h,v)=g \cdot \phi(h,v)$ for all $g \in \GSpint(V)$ and $(h,v) \in X$ and so this $\phi$ is a continuous $\GSpint(V)$-equivariant Nash map. Thus  by the Frobenius descent we have
\[
\Scal^*(X)^{\GSpint(V), \chi}\simeq \Scal^*(\GSpin(V)\times\{e\})^{\GSpint(V)_e, \chi},
\]
where the left-hand side is a subspace of $\Scal^*(\GSpin(V)\times V)^{\GSpint(V), \chi}$. But clearly
\[
\Scal^*(\GSpin(V)\times\{e\})^{\GSpint(V)_e, \chi}\simeq \Scal'(\GSpint(V))^{\GSpint(W), \chi}.
\]
The proposition follows.
\end{proof}

\subsection{Reducing to classical group situation}
It is now enough to show
\[
\Scal^*(\GSpin(V)\times V)^{\GSpint(V), \chi}=0.
\]
Arguing in the same way as the $\GPin$ case, this vanishing assertion reduces to
\[
\Scal^*(g\,\Ucal_g\times V)^{\GSpint(V)_g, \chi}=0
\]
for all semisimple $g\in\GSpin(V)$. Since the canonical projection $P$ is bijective on $g\,\Ucal_g$, we have the natural isomorphism
\[
\Scal^*(g\,\Ucal_g\times V)^{\GSpint(V)_g, \chi}\simeq\Scal'(P(g\,\Ucal_g)\times V)^{P(\GSpint(V)_g), \chi}.
\]
Since we have
\[
\Scal^*(P(g\,\Ucal_g)\times V)^{P(\GSpint(V)_g), \chi}
\subseteq\Scal'(P(\GSpin(V)_g)\times V)^{P(\GSpint(V)_g), \chi},
\]
it is enough to show
\[
\Scal'(P(\GSpin(V)_g)\times V)^{P(\GSpint(V)_g), \chi}=0.
\]
We know $P(\GSpin(V)_g)$ is as in Lemma \ref{lemma:centralizer_GSpin} and $P(\GSpint(V)_g)$ is generated by $P(\GSpin(V)_g)$ and the element $\gamma\beta$, where $\gamma=(\gamma_1,\dots,\gamma_m,\gamma_+,\gamma_-)\in\OO(V_1)\times\cdots\times\OO(V_m)\times\OO(V_+)\times\OO(V_-)
$
such that $\gamma h^{-1}\gamma^{-1}=h$, where  $h:=P(g)\in\SO(V)$ as before (see also \cite[Equation 8.1]{ET2023}. Note that since the orthogonal factor of $P(\GSpin(V)_g)$ is $\SO(V_+)\times\SO(V_-)$, we always choose $\gamma_{\pm}=r_{e_{\pm}}^{k_{\pm}}$. The remaining part of the proof is identical to the proof in the GPin case, and so the proof is complete.

In conclusion, we have shown the vanishing assertion of \eqref{eq: vanishing} which gives Theorem \ref{thm:main} for the groups $\GSpin (V)$ and $\GSpin (W)$.

\bibliographystyle{alpha}
\bibliography{MultOne}
\end{document}